\documentclass[final,notitlepage,onecolumn]{article}%
\usepackage{natbib}
\usepackage{amsmath,amsthm}
\usepackage{graphicx}
\usepackage{amsfonts}
\usepackage{amssymb}
\usepackage{amsmath}%
\providecommand{\U}[1]{\protect\rule{.1in}{.1in}}
\providecommand{\U}[1]{\protect\rule{.1in}{.1in}}

\newtheorem{theorem}{Theorem}
\theoremstyle{remark}

\begin{document}

\begin{center}
\vspace{0.5in}%


\vspace{0.25in}

{\Huge DEPARTMENT OF MATHEMATICAL SCIENCES}

\smallskip

{\Huge TECHNICAL REPORT SERIES}

\vspace{0.5in}

{\LARGE The Generalized Law of Total Covariance}

\bigskip

by

\bigskip

Charles W. Champ and Andrew V. Sills

Department of Mathematical Sciences

Georgia Southern University, Statesboro, GA 30460-8093

\vspace{2in}

Number 2022-001

Submitted: May 31, 2022

\copyright \ 2022
\end{center}

\pagebreak

\section{\textbf{Introduction}}

The covariance $Cov\left(  Y_{1},Y_{2}\right)  $\ between the random
measurements $Y_{1}$ and $Y_{2}$ to be taken on each individual in the
population is the average value of the random variable $\left[  Y_{1}-E\left(
Y_{1}\right)  \right]  \left[  Y_{2}-E\left(  Y_{2}\right)  \right]  $ over
all individuals in the population. That is,%
\begin{align*}
Cov\left(  Y_{1},Y_{2}\right)   &  =Cov_{Y_{1},Y_{2}}\left(  Y_{1}%
,Y_{2}\right)  =E_{Y_{1},Y_{2}}\left(  \left[  Y_{1}-E\left(  Y_{1}\right)
\right]  \left[  Y_{2}-E\left(  Y_{2}\right)  \right]  \right) \\
&  =\int\nolimits_{-\infty}^{\infty}\int\nolimits_{-\infty}^{\infty}\left[
y_{1}-E\left(  Y_{1}\right)  \right]  \left[  y_{2}-E\left(  Y_{2}\right)
\right]  f_{Y_{1},Y_{2}}\left(  y_{1},y_{2}\right)  dy_{1}dy_{2}%
\end{align*}
if both $Y_{1}$ and $Y_{2}$ are continuous random variables. The subscript
\textquotedblleft$Y_{1},Y_{2}$\textquotedblright\ is a more explicit way of
indicating that the covariance is a parameter of the joint distribution of
$Y_{1}$ and $Y_{2}$. If $Y_{1}$ and/or $Y_{2}$ is a discrete random variable,
integration is replaced with summation. For convenience, we let $\mathbf{Y}%
=\left[  Y_{1},Y_{2}\right]  ^{\mathbf{T}}$. We will also express the
$Cov\left(  Y_{1},Y_{2}\right)  $ by $Cov\left(  \mathbf{Y}\right)  $ and
$Cov_{\mathbf{Y}}\left(  \mathbf{Y}\right)  $. It is not difficult to show
that%
\[
Cov_{\mathbf{Y}}\left(  \mathbf{Y}\right)  =E_{\mathbf{Y}}\left(  Y_{1}%
Y_{2}\right)  -E_{\mathbf{Y}}\left(  Y_{1}\right)  E_{\mathbf{Y}}\left(
Y_{2}\right)  \text{.}%
\]
Supposing that the random vector $\mathbf{Y}$ and the random variable $X_{1}$
have a joint distribution, the Total Law of Covariance states that%
\[
Cov_{\mathbf{Y}}\left(  \mathbf{Y}\right)  =E_{X_{1}}\left[  Cov_{\mathbf{Y}%
\left\vert X_{1}\right.  }\left(  Y_{1},Y_{2}\right)  \right]  +Cov_{X_{1}%
}\left(  E_{\mathbf{Y}\left\vert X_{1}\right.  }\left(  Y_{1}\right)
,E_{\mathbf{Y}\left\vert X_{1}\right.  }\left(  Y_{2}\right)  \right)
\text{.}%
\]
The total law of covariance is also referred to as covariance decomposition
formula or conditional covariance formula.

For the case in which $Y_{1}=Y_{2}=Y$, then $Cov\left(  Y,Y\right)  =V\left(
Y\right)  $. For this case it follows that%
\begin{equation}
V\left(  Y\right)  =E_{X_{1}}\left[  V_{Y\left\vert X_{1}\right.  }\left(
Y_{1},Y_{2}\right)  \right]  +V_{X_{1}}\left[  E_{Y\left\vert X_{1}\right.
}\left(  Y\right)  \right]  \text{.} \label{LTV}%
\end{equation}
This special case is known as the total law of variance. The law of total
variance is also referred to as the \textit{variance decomposition formula},
the \textit{conditional variance formulas}, or \textit{Law of Iterated
Variances} -- known as \textit{Eve's Law} (see~\citet{BH}).
~\citet{BE},~\citet{R},~\cite{WMS}. The preceding references provide proofs
of~\eqref{LTV}. Applications of~\eqref{LTV} in actuarial science can be found
in~\cite{MD}.

It is our interest to examine the unconditional mean $E_{\mathbf{Y}}\left[
\left(  \mathbf{Y}\right)  \right]  $ and the unconditional variance
$Cov_{\mathbf{Y}}\left(  \mathbf{Y}\right)  $ of the distribution of
$\mathbf{Y}$ with respect to the joint distribution of $\mathbf{Y}$ and the
$k$ random variables $X_{1},\ldots,X_{k}$. In the next section, we will
present results obtained by Champ and Sills (2019) for the case in which
$Y_{1}=Y_{2}=Y$. They refer to their results as the \emph{generalized law of
total variance}. In Section 3, we present a generalization of the law of total covariance.

\smallskip

\noindent Key Words: Law of total variance, conditional variance and
covariance, conditional expectation, Eve's law, variance decomposition,
covariance decompostion

\section{The Generalized Law of Total Expectation and Variance}

Champ and Sills (2019) present and prove the following theorems which we
include for completeness in this artile. The first theorem is about
expectation. This is a generalization of the what is known as the law of total
expectation, the law of iterated expectations, the tower rule, Adam's law, and
the smoothing theorem. Champ and Sills (2019) refer to their results as the
generalized law of total expection.

\smallskip

\begin{theorem}
Assume that $Y$ and the $k$ random variables $X_{1},\ldots,X_{k}$ have a joint
distribution. Then%
\[
E\left(  Y\right)  =E_{X_{1}}\left[  E_{X_{2}\left\vert \mathbf{X}_{1}\right.
}\cdots\left[  E_{X_{k}\left\vert \mathbf{X}_{k-1}\right.  }\left[
E_{Y\left\vert \mathbf{X}_{k}\right.  }\left(  Y\right)  \right]  \right]
\cdots\right]  \text{,}%
\]
where $\mathbf{X}_{i}=\left[  X_{1},\ldots,X_{i}\right]  ^{\mathbf{T}}$ for
$i=1,\ldots,k$.
\end{theorem}

\smallskip

Champ and Sills (2019) refer to the next theorem as the \textit{Generalized
Law of Total Variance}.

\smallskip

\begin{theorem}
\label{GLTV} Assume that $Y$ and the $k$ random variables $X_{1},\ldots,X_{k}
$ have a joint distribution. Then%
\begin{equation}
\label{GLTVeq}V\left(  Y\right)  =\sum\nolimits_{i=1}^{k+1}Q_{X_{1}%
}^{\left\langle i,1\right\rangle }\left[  Q_{X_{2}\left\vert X_{1}\right.
}^{\left\langle i,2\right\rangle }\left[  \cdots\left[  Q_{Y\left\vert
X_{1},\ldots,X_{k}\right.  }^{\left\langle i,k+1\right\rangle }\left(
Y\right)  \right]  \cdots\right]  \right]  \text{,}%
\end{equation}
where the operator $Q^{\left\langle i,j\right\rangle }$ is the variance
operator $V$ if $i=j$ and the expectation operator $E$ if $i\neq j$, for
$i,j=1,\ldots,k+1$.
\end{theorem}

\smallskip

For $k=1$, we have%
\[
V\left(  Y\right)  =E_{X_{1}}\left[  V_{Y\left\vert \mathbf{X}_{1}\right.
}\left(  Y\right)  \right]  +V_{X_{1}}\left[  E_{Y\left\vert \mathbf{X}%
_{1}\right.  }\left(  Y\right)  \right]  \text{.}%
\]
This is the law of total variance. For the case in which $k=2$, the variance
of $Y$ can be expressed as%
\[
V\left(  Y\right)  =E_{X_{1}}\left[  E_{X_{2}\left\vert \mathbf{X}_{1}\right.
}\left[  V_{Y\left\vert \mathbf{X}_{2}\right.  }\left(  Y\right)  \right]
\right]  +E_{X_{1}}\left[  V_{X_{2}\left\vert \mathbf{X}_{1}\right.  }\left[
E_{Y\left\vert \mathbf{X}_{2}\right.  }\left(  Y\right)  \right]  \right]
+V_{X_{1}}\left[  E_{X_{2}\left\vert \mathbf{X}_{1}\right.  }\left[
E_{Y\left\vert \mathbf{X}_{2}\right.  }\left(  Y\right)  \right]  \right]
\text{.}%
\]
An equivalent result is given in \citet{GA}. \citet{KPW} provide an example
applying the law of total variance. Also, \citet{M} provides an example of
that applies the law of total variance.

\smallskip

\section{The Generalized Law of Total Covariance}

\begin{theorem}
If $p\times1$ random vector $\mathbf{Y}$ and the $k$ random variables
$X_{1},\ldots,X_{k}$ have a joint distribution, then the expection of
$u\left(  \mathbf{Y}\right)  $, a real valued function of $\mathbf{Y}$, can be
expressed as%
\[
E\left[  u\left(  \mathbf{Y}\right)  \right]  =E_{X_{1}}\left[  E_{X_{2}%
\left\vert \mathbf{X}_{1}\right.  }\cdots\left[  E_{X_{k}\left\vert
\mathbf{X}_{k-1}\right.  }\left[  E_{\mathbf{Y}\left\vert \mathbf{X}%
_{k}\right.  }\left[  u\left(  \mathbf{Y}\right)  \right]  \right]  \right]
\cdots\right]  \text{,}%
\]
where $\mathbf{X}_{i}=\left[  X_{1},\ldots,X_{i}\right]  ^{\mathbf{T}}$ for
$i=1,\ldots,k$.
\end{theorem}

\noindent Proof: We see for $k=1$ that%
\begin{align*}
E\left[  u\left(  \mathbf{Y}\right)  \right]   &  =\int\nolimits_{-\infty
}^{\infty}\int\nolimits_{-\infty}^{\infty}u\left(  \mathbf{y}\right)
f_{\mathbf{Y}}\left(  \mathbf{y}\right)  d\mathbf{y}=\int\nolimits_{-\infty
}^{\infty}\int\nolimits_{-\infty}^{\infty}u\left(  \mathbf{y}\right)  \left[
\int\nolimits_{-\infty}^{\infty}f_{\mathbf{Y}}\left(  \mathbf{y},x_{1}\right)
dx_{1}\right]  d\mathbf{y}\\
&  =\int\nolimits_{-\infty}^{\infty}\int\nolimits_{-\infty}^{\infty}u\left(
\mathbf{y}\right)  \left[  \int\nolimits_{-\infty}^{\infty}f_{\mathbf{Y}%
\left\vert X_{1}\right.  }\left(  \mathbf{y}\left\vert x_{1}\right.  \right)
f_{X_{1}}\left(  x_{1}\right)  dx_{1}\right]  d\mathbf{y}\\
&  =\int\nolimits_{-\infty}^{\infty}\left[  \int\nolimits_{-\infty}^{\infty
}\int\nolimits_{-\infty}^{\infty}u\left(  \mathbf{y}\right)  f_{\mathbf{Y}%
\left\vert X_{1}\right.  }\left(  \mathbf{y}\left\vert x_{1}\right.  \right)
d\mathbf{y}\right]  f_{X_{1}}\left(  x_{1}\right)  dx_{1}\\
&  =\int\nolimits_{-\infty}^{\infty}E_{\mathbf{Y}\left\vert x_{1}\right.
}\left[  u\left(  \mathbf{Y}\right)  \right]  f_{X_{1}}\left(  x_{1}\right)
dx_{1}\\
&  =E_{X_{1}}\left[  E_{\mathbf{Y}\left\vert \mathbf{X}_{1}\right.  }\left[
u\left(  \mathbf{Y}\right)  \right]  \right]  \text{.}%
\end{align*}
Thus, the theorem is true for $k=1$.

Now assume the theorem is true for $k>1$. Further, assume that $\mathbf{Y}$
and the $k+1$ random variables $X_{1},\ldots,X_{k+1}$ have a joint
distribution. Observe that%
\begin{align*}
E_{\mathbf{Y}\left\vert \mathbf{X}_{k}\right.  }\left[  u\left(
\mathbf{Y}\right)  \right]   &  =\int\nolimits_{-\infty}^{\infty}%
\int\nolimits_{-\infty}^{\infty}u\left(  \mathbf{y}\right)  f_{\mathbf{Y}%
\left\vert \mathbf{X}_{k}\right.  }\left(  \mathbf{y}\left\vert \mathbf{x}%
_{k}\right.  \right)  d\mathbf{y}=\int\nolimits_{-\infty}^{\infty}%
\int\nolimits_{-\infty}^{\infty}u\left(  \mathbf{y}\right)  \left[
\frac{f_{\mathbf{Y},\mathbf{X}_{k}}\left(  \mathbf{y},\mathbf{x}_{k}\right)
}{f_{\mathbf{X}_{k}}\left(  \mathbf{x}_{k}\right)  }\right]  d\mathbf{y}\\
&  =\int\nolimits_{-\infty}^{\infty}\int\nolimits_{-\infty}^{\infty}u\left(
\mathbf{y}\right)  \left[  \frac{\int\nolimits_{-\infty}^{\infty}%
f_{\mathbf{Y},\mathbf{X}_{k},X_{k+1}}\left(  \mathbf{y},\mathbf{x}_{k}%
,x_{k+1}\right)  dx_{k+1}}{f_{\mathbf{X}_{k}}\left(  \mathbf{x}_{k}\right)
}\right]  d\mathbf{y}\\
&  =\int\nolimits_{-\infty}^{\infty}\int\nolimits_{-\infty}^{\infty}u\left(
\mathbf{y}\right)  \left[  \frac{\int\nolimits_{-\infty}^{\infty}%
f_{\mathbf{Y}\left\vert \mathbf{X}_{k},X_{k+1}\right.  }\left(  \mathbf{y}%
\left\vert \mathbf{x}_{k},x_{k+1}\right.  \right)  f_{X_{k+1}\left\vert
\mathbf{X}_{k}\right.  }\left(  x_{k+1}\left\vert \mathbf{x}_{k}\right.
\right)  f_{\mathbf{X}_{k}}\left(  \mathbf{x}_{k}\right)  dx_{k+1}%
}{f_{\mathbf{X}_{k}}\left(  \mathbf{x}_{k}\right)  }\right]  d\mathbf{y}\\
&  =\int\nolimits_{-\infty}^{\infty}\int\nolimits_{-\infty}^{\infty}u\left(
\mathbf{y}\right)  \left[  \int\nolimits_{-\infty}^{\infty}f_{\mathbf{Y}%
\left\vert \mathbf{X}_{k},X_{k+1}\right.  }\left(  \mathbf{y}\left\vert
\mathbf{x}_{k},x_{k+1}\right.  \right)  f_{X_{k+1}\left\vert \mathbf{X}%
_{k}\right.  }\left(  x_{k+1}\left\vert \mathbf{x}_{k}\right.  \right)
dx_{k+1}\right]  d\mathbf{y}\\
&  =\int\nolimits_{-\infty}^{\infty}\left[  \int\nolimits_{-\infty}^{\infty
}\int\nolimits_{-\infty}^{\infty}u\left(  \mathbf{y}\right)  f_{\mathbf{Y}%
\left\vert \mathbf{X}_{k},X_{k+1}\right.  }\left(  \mathbf{y}\left\vert
\mathbf{x}_{k},x_{k+1}\right.  \right)  d\mathbf{y}\right]  f_{X_{k+1}%
\left\vert \mathbf{X}_{k}\right.  }\left(  x_{k+1}\left\vert \mathbf{x}%
_{k}\right.  \right)  dx_{k+1}\\
&  =\int\nolimits_{-\infty}^{\infty}E_{\mathbf{Y}\left\vert \mathbf{X}%
_{k},x_{k+1}\right.  }\left[  u\left(  \mathbf{Y}\right)  \right]
f_{X_{k+1}\left\vert \mathbf{X}_{k}\right.  }\left(  x_{k+1}\left\vert
\mathbf{x}_{k}\right.  \right)  dx_{k+1}\\
&  =E_{X_{k+1}\left\vert \mathbf{X}_{k}\right.  }\left[  E_{\mathbf{Y}%
\left\vert \mathbf{X}_{k},X_{k+1}\right.  }\left[  u\left(  \mathbf{Y}\right)
\right]  \right] \\
&  =E_{X_{k+1}\left\vert \mathbf{X}_{k}\right.  }\left[  E_{\mathbf{Y}%
\left\vert \mathbf{X}_{k+1}\right.  }\left[  u\left(  \mathbf{Y}\right)
\right]  \right]  \text{.}%
\end{align*}
Substituting $E_{\mathbf{Y}\left\vert \mathbf{X}_{k}\right.  }\left[  u\left(
\mathbf{Y}\right)  \right]  $ with $E_{X_{k+1}\left\vert \mathbf{X}%
_{k}\right.  }\left[  E_{\mathbf{Y}\left\vert \mathbf{X}_{k+1}\right.
}\left[  u\left(  \mathbf{Y}\right)  \right]  \right]  $ in the theorem
results in%
\[
E\left[  u\left(  \mathbf{Y}\right)  \right]  =E_{X_{1}}\left[  E_{X_{2}%
\left\vert \mathbf{X}_{1}\right.  }\cdots\left[  E_{X_{k}\left\vert
\mathbf{X}_{k-1}\right.  }\left[  E_{X_{k+1}\left\vert \mathbf{X}_{k}\right.
}\left[  E_{\mathbf{Y}\left\vert \mathbf{X}_{k+1}\right.  }\left[  u\left(
\mathbf{Y}\right)  \right]  \right]  \right]  \right]  \cdots\right]  \text{.}%
\]
Hence, the theorem is true for $k+1$. Therefore, by the Axiom of Induction,
the theorem is true for all positive integers $k$.$\blacksquare$

\smallskip

\begin{theorem}
(Generalized Law of Total Covariance) Assume that $\mathbf{Y}=\left[
Y_{1},Y_{2}\right]  ^{\mathbf{T}}$ and the $k$ random variables $X_{1}%
,\ldots,X_{k}$ have a joint distribution. Then%
\begin{align*}
Cov\left(  \mathbf{Y}\right)   &  =Cov\left(  Y_{1},Y_{2}\right) \\
&  =Cov_{X_{1}}\left(  E_{X_{2}\left\vert \mathbf{X}_{1}\right.  }\left[
\cdots\left[  E_{\mathbf{Y}\left\vert \mathbf{X}_{k}\right.  }\left(
Y_{1}\right)  \right]  \cdots\right]  ,E_{X_{2}\left\vert \mathbf{X}%
_{1}\right.  }\left[  \cdots\left[  E_{\mathbf{Y}\left\vert \mathbf{X}%
_{k}\right.  }\left(  Y_{2}\right)  \right]  \cdots\right]  \right) \\
&  +E_{X_{1}}\left[  Cov_{X_{2}\left\vert \mathbf{X}_{1}\right.  }\left(
E_{X_{3}\left\vert \mathbf{X}_{2}\right.  }\left[  \cdots\left[
E_{\mathbf{Y}\left\vert \mathbf{X}_{k}\right.  }\left(  Y_{1}\right)  \right]
\cdots\right]  ,E_{X_{3}\left\vert \mathbf{X}_{2}\right.  }\left[
\cdots\left[  E_{\mathbf{Y}\left\vert \mathbf{X}_{k}\right.  }\left(
Y_{2}\right)  \right]  \cdots\right]  \right)  \right] \\
&  \vdots\\
&  +E_{X_{1}}\left[  E_{X_{2}\left\vert \mathbf{X}_{1}\right.  }\left[
E_{X_{3}\left\vert \mathbf{X}_{2}\right.  }\left[  \cdots\left[
Cov_{\mathbf{Y}\left\vert \mathbf{X}_{k}\right.  }\left(  Y_{1},Y_{2}\right)
\right]  \cdots\right]  \right]  \right]  \text{,}%
\end{align*}
where $\mathbf{X}_{i}=\left[  X_{1},\ldots,X_{i}\right]  ^{\mathbf{T}}$ for
$i=1,\ldots,k$.
\end{theorem}

\noindent Proof: As previously stated,%
\[
Cov_{\mathbf{Y}}\left(  \mathbf{Y}\right)  =E_{\mathbf{Y}}\left(  Y_{1}%
Y_{2}\right)  -E_{\mathbf{Y}}\left(  Y_{1}\right)  E_{\mathbf{Y}}\left(
Y_{2}\right)  \text{.}%
\]
For $k=1$, using theorem 3 with $u\left(  \mathbf{Y}\right)  =Y_{1}Y_{2}$, we
can write%
\[
E_{\mathbf{Y}}\left(  Y_{1}Y_{2}\right)  =E_{X_{1}}\left[  \left[
E_{\mathbf{Y}\left\vert X_{1}\right.  }\left(  Y_{1}Y_{2}\right)  \right]
\right]
\]
and with $u\left(  \mathbf{Y}\right)  =Y_{i}$ for $i=1,2$,%
\[
E_{\mathbf{Y}}\left(  Y_{1}\right)  E_{\mathbf{Y}}\left(  Y_{2}\right)
=E_{X_{1}}\left[  \left[  E_{\mathbf{Y}\left\vert X_{1}\right.  }\left(
Y_{1}\right)  \right]  \right]  E_{X_{1}}\left[  \left[  E_{\mathbf{Y}%
\left\vert X_{1}\right.  }\left(  Y_{2}\right)  \right]  \right]  \text{.}%
\]
Thus,%
\begin{align*}
Cov_{\mathbf{Y}}\left(  \mathbf{Y}\right)   &  =E_{X_{1}}\left[  \left[
E_{\mathbf{Y}\left\vert X_{1}\right.  }\left(  Y_{1}Y_{2}\right)  \right]
\right]  -E_{X_{1}}\left[  \left[  E_{\mathbf{Y}\left\vert X_{1}\right.
}\left(  Y_{1}\right)  \right]  \right]  E_{X_{1}}\left[  \left[
E_{\mathbf{Y}\left\vert X_{1}\right.  }\left(  Y_{2}\right)  \right]  \right]
\\
&  =E_{X_{1}}\left[  \left[  E_{\mathbf{Y}\left\vert X_{1}\right.  }\left(
Y_{1}Y_{2}\right)  \right]  -\left[  E_{\mathbf{Y}\left\vert X_{1}\right.
}\left(  Y_{1}\right)  \right]  \left[  E_{\mathbf{Y}\left\vert X_{1}\right.
}\left(  Y_{2}\right)  \right]  \right]  \\
&  +E_{X_{1}}\left[  \left[  E_{\mathbf{Y}\left\vert X_{1}\right.  }\left(
Y_{1}\right)  \right]  \left[  E_{\mathbf{Y}\left\vert X_{1}\right.  }\left(
Y_{2}\right)  \right]  \right]  -E_{X_{1}}\left[  \left[  E_{\mathbf{Y}%
\left\vert X_{1}\right.  }\left(  Y_{1}\right)  \right]  \right]  E_{X_{1}%
}\left[  \left[  E_{\mathbf{Y}\left\vert X_{1}\right.  }\left(  Y_{2}\right)
\right]  \right]  \\
&  =E_{X_{1}}\left[  Cov\left(  E_{\mathbf{Y}\left\vert X_{1}\right.  }\left(
Y_{1}\right)  ,E_{\mathbf{Y}\left\vert X_{1}\right.  }\left(  Y_{2}\right)
\right)  \right]  +Cov_{X_{1}}\left(  E_{\mathbf{Y}\left\vert X_{1}\right.
}\left(  Y_{1}\right)  ,E_{\mathbf{Y}\left\vert X_{1}\right.  }\left(
Y_{2}\right)  \right)  \text{.}%
\end{align*}
Hence, the theorem is true forrem is true for $k>1$. Further, assume that
$\mathbf{Y}$ and the $k+1$ random variables $X_{1},\ldots,X_{k},X_{k+1}$ have
a joint distribution. Observe that%
\[
E_{\mathbf{Y}\left\vert \mathbf{X}_{k}\right.  }\left(  Y_{i}\right)
=E_{X_{k+1}\left\vert \mathbf{X}_{k}\right.  }\left[  E_{\mathbf{Y}\left\vert
\mathbf{X}_{k+1}\right.  }\left(  Y_{i}\right)  \right]
\]
for $i=1,2$ which follows from the proof of theorem 3. Again using the results
of theorem 3, we have that%
\begin{align*}
Cov_{\mathbf{Y}\left\vert \mathbf{X}_{k}\right.  }\left(  Y_{1},Y_{2}\right)
&  =E_{\mathbf{Y}\left\vert \mathbf{X}_{k}\right.  }\left(  Y_{1}Y_{2}\right)
-E_{\mathbf{Y}\left\vert \mathbf{X}_{k}\right.  }\left(  Y_{1}\right)
E_{\mathbf{Y}\left\vert \mathbf{X}_{k}\right.  }\left(  Y_{2}\right)  \\
&  =E_{X_{k+1}\left\vert \mathbf{X}_{k}\right.  }\left[  E_{\mathbf{Y}%
\left\vert \mathbf{X}_{k+1}\right.  }\left(  Y_{1}Y_{2}\right)  -E_{\mathbf{Y}%
\left\vert \mathbf{X}_{k+1}\right.  }\left(  Y_{1}\right)  E_{\mathbf{Y}%
\left\vert \mathbf{X}_{k+1}\right.  }\left(  Y_{2}\right)  \right]  \\
&  +E_{X_{k+1}\left\vert \mathbf{X}_{k}\right.  }\left[  E_{\mathbf{Y}%
\left\vert \mathbf{X}_{k+1}\right.  }\left(  Y_{1}\right)  E_{\mathbf{Y}%
\left\vert \mathbf{X}_{k+1}\right.  }\left(  Y_{2}\right)  \right]  \\
&  -E_{X_{k+1}\left\vert \mathbf{X}_{k}\right.  }\left[  E_{\mathbf{Y}%
\left\vert \mathbf{X}_{k+1}\right.  }\left(  Y_{1}\right)  \right]
E_{X_{k+1}\left\vert \mathbf{X}_{k}\right.  }\left[  E_{\mathbf{Y}\left\vert
\mathbf{X}_{k+1}\right.  }\left(  Y_{2}\right)  \right]  \\
&  =E_{X_{k+1}\left\vert \mathbf{X}_{k}\right.  }\left[  Cov_{\mathbf{Y}%
\left\vert \mathbf{X}_{k+1}\right.  }\left(  Y_{1},Y_{2}\right)  \right]  \\
&  +E_{X_{k+1}\left\vert \mathbf{X}_{k}\right.  }\left[  Cov_{\mathbf{Y}%
\left\vert \mathbf{X}_{k+1}\right.  }\left(  E_{\mathbf{Y}\left\vert
\mathbf{X}_{k+1}\right.  }\left(  Y_{1}\right)  ,E_{\mathbf{Y}\left\vert
\mathbf{X}_{k+1}\right.  }\left(  Y_{2}\right)  \right)  \right]  \\
&  =E_{X_{k+1}\left\vert \mathbf{X}_{k}\right.  }\left[  Cov_{\mathbf{Y}%
\left\vert \mathbf{X}_{k+1}\right.  }\left(  E_{\mathbf{Y}\left\vert
\mathbf{X}_{k+1}\right.  }\left(  Y_{1}\right)  ,E_{\mathbf{Y}\left\vert
\mathbf{X}_{k+1}\right.  }\left(  Y_{2}\right)  \right)  \right]  \\
&  +E_{X_{k+1}\left\vert \mathbf{X}_{k}\right.  }\left[  Cov_{\mathbf{Y}%
\left\vert \mathbf{X}_{k+1}\right.  }\left(  Y_{1},Y_{2}\right)  \right]
\text{.}%
\end{align*}
Substituting $E_{\mathbf{Y}\left\vert \mathbf{X}_{k}\right.  }\left(
Y_{i}\right)  $ with $E_{X_{k+1}\left\vert \mathbf{X}_{k}\right.  }\left[
E_{\mathbf{Y}\left\vert \mathbf{X}_{k+1}\right.  }\left(  Y_{i}\right)
\right]  $ for $i=1,2$ and $Cov_{\mathbf{Y}\left\vert \mathbf{X}_{k}\right.
}\left(  Y_{1},Y_{2}\right)  $ with%
\[
E_{X_{k+1}\left\vert \mathbf{X}_{k}\right.  }\left[  Cov_{\mathbf{Y}\left\vert
\mathbf{X}_{k+1}\right.  }\left(  E_{\mathbf{Y}\left\vert \mathbf{X}%
_{k+1}\right.  }\left(  Y_{1}\right)  ,E_{\mathbf{Y}\left\vert \mathbf{X}%
_{k+1}\right.  }\left(  Y_{2}\right)  \right)  \right]  +E_{X_{k+1}\left\vert
\mathbf{X}_{k}\right.  }\left[  Cov_{\mathbf{Y}\left\vert \mathbf{X}%
_{k+1}\right.  }\left(  Y_{1},Y_{2}\right)  \right]
\]
yields%
\begin{align*}
&  Cov\left(  \mathbf{Y}\right)  \\
&  =Cov_{X_{1}}\left(  E_{X_{2}\left\vert \mathbf{X}_{1}\right.  }\left[
\cdots\left[  E_{X_{k+1}\left\vert \mathbf{X}_{k}\right.  }\left[
E_{\mathbf{Y}\left\vert \mathbf{X}_{k+1}\right.  }\left(  Y_{1}\right)
\right]  \right]  \cdots\right]  ,E_{X_{2}\left\vert \mathbf{X}_{1}\right.
}\left[  \cdots\left[  E_{X_{k+1}\left\vert \mathbf{X}_{k}\right.  }\left[
E_{\mathbf{Y}\left\vert \mathbf{X}_{k+1}\right.  }\left(  Y_{2}\right)
\right]  \right]  \cdots\right]  \right)  \\
&  +E_{X_{1}}\left[  Cov_{X_{2}\left\vert \mathbf{X}_{1}\right.  }\left(
\begin{array}
[c]{c}%
E_{X_{3}\left\vert \mathbf{X}_{2}\right.  }\left[  \cdots\left[
E_{X_{k+1}\left\vert \mathbf{X}_{k}\right.  }\left[  E_{\mathbf{Y}\left\vert
\mathbf{X}_{k+1}\right.  }\left(  Y_{1}\right)  \right]  \right]
\cdots\right]  \\
,E_{X_{3}\left\vert \mathbf{X}_{2}\right.  }\left[  \cdots\left[
E_{X_{k+1}\left\vert \mathbf{X}_{k}\right.  }\left[  E_{\mathbf{Y}\left\vert
\mathbf{X}_{k+1}\right.  }\left(  Y_{2}\right)  \right]  \right]
\cdots\right]
\end{array}
\right)  \right]  \\
&  \vdots\\
&  +E_{X_{1}}\left[  E_{X_{2}\left\vert \mathbf{X}_{1}\right.  }\left[
E_{X_{3}\left\vert \mathbf{X}_{2}\right.  }\left[  \cdots\left[
\begin{array}
[c]{c}%
E_{X_{k+1}\left\vert \mathbf{X}_{k}\right.  }\left[  Cov_{\mathbf{Y}\left\vert
\mathbf{X}_{k+1}\right.  }\left(
\begin{array}
[c]{c}%
E_{\mathbf{Y}\left\vert \mathbf{X}_{k+1}\right.  }\left(  Y_{1}\right)  ,\\
E_{\mathbf{Y}\left\vert \mathbf{X}_{k+1}\right.  }\left(  Y_{2}\right)
\end{array}
\right)  \right]  \\
+E_{X_{k+1}\left\vert \mathbf{X}_{k}\right.  }\left[  Cov_{\mathbf{Y}%
\left\vert \mathbf{X}_{k+1}\right.  }\left(  Y_{1},Y_{2}\right)  \right]
\end{array}
\right]  \cdots\right]  \right]  \right]  \\
&  =Cov_{X_{1}}\left(
\begin{array}
[c]{c}%
E_{X_{2}\left\vert \mathbf{X}_{1}\right.  }\left[  \cdots\left[
E_{X_{k+1}\left\vert \mathbf{X}_{k}\right.  }\left[  E_{\mathbf{Y}\left\vert
\mathbf{X}_{k+1}\right.  }\left(  Y_{1}\right)  \right]  \right]
\cdots\right]  \\
,E_{X_{2}\left\vert \mathbf{X}_{1}\right.  }\left[  \cdots\left[
E_{X_{k+1}\left\vert \mathbf{X}_{k}\right.  }\left[  E_{\mathbf{Y}\left\vert
\mathbf{X}_{k+1}\right.  }\left(  Y_{2}\right)  \right]  \right]
\cdots\right]
\end{array}
\right)  \\
&  +E_{X_{1}}\left[  Cov_{X_{2}\left\vert \mathbf{X}_{1}\right.  }\left(
\begin{array}
[c]{c}%
E_{X_{3}\left\vert \mathbf{X}_{2}\right.  }\left[  \cdots\left[
E_{X_{k+1}\left\vert \mathbf{X}_{k}\right.  }\left[  E_{\mathbf{Y}\left\vert
\mathbf{X}_{k+1}\right.  }\left(  Y_{1}\right)  \right]  \right]
\cdots\right]  \\
,E_{X_{3}\left\vert \mathbf{X}_{2}\right.  }\left[  \cdots\left[
E_{X_{k+1}\left\vert \mathbf{X}_{k}\right.  }\left[  E_{\mathbf{Y}\left\vert
\mathbf{X}_{k+1}\right.  }\left(  Y_{2}\right)  \right]  \right]
\cdots\right]
\end{array}
\right)  \right]  \\
&  \vdots\\
&  +E_{X_{1}}\left[  E_{X_{2}\left\vert \mathbf{X}_{1}\right.  }\left[
E_{X_{3}\left\vert \mathbf{X}_{2}\right.  }\left[  \cdots\left[
E_{X_{k+1}\left\vert \mathbf{X}_{k}\right.  }\left[  Cov_{\mathbf{Y}\left\vert
\mathbf{X}_{k+1}\right.  }\left(
\begin{array}
[c]{c}%
E_{\mathbf{Y}\left\vert \mathbf{X}_{k+1}\right.  }\left(  Y_{1}\right)  \\
,E_{\mathbf{Y}\left\vert \mathbf{X}_{k+1}\right.  }\left(  Y_{2}\right)
\end{array}
\right)  \right]  \right]  \cdots\right]  \right]  \right]  \\
&  +E_{X_{1}}\left[  E_{X_{2}\left\vert \mathbf{X}_{1}\right.  }\left[
E_{X_{3}\left\vert \mathbf{X}_{2}\right.  }\left[  \cdots\left[
E_{X_{k+1}\left\vert \mathbf{X}_{k}\right.  }E_{X_{k+1}\left\vert
\mathbf{X}_{k}\right.  }\left[  Cov_{\mathbf{Y}\left\vert \mathbf{X}%
_{k+1}\right.  }\left(  Y_{1},Y_{2}\right)  \right]  \right]  \cdots\right]
\right]  \right]  \text{.}%
\end{align*}
Hence, the theorem is true for $k+1$. Therefore, by the Axiom of Induction,
the theorem is true for all positive integers $k$.$\blacksquare$

\smallskip

Klugman, Panjer, and Willmot (2019) and Maheshwari (2020) provide examples
applying the law of total covariance.

\smallskip

\newpage

\begin{theorem}
The generalized law of total variance is a special case of the generalized law
of total covariance.
\end{theorem}

\noindent Proof: Suppose that $Y_{1}=Y_{2}=Y$. We have%
\begin{align*}
V\left(  Y\right)   &  =Cov\left(  Y,Y\right) \\
&  =Cov_{X_{1}}\left(  E_{X_{2}\left\vert \mathbf{X}_{1}\right.  }\left[
\cdots\left[  E_{\mathbf{Y}\left\vert \mathbf{X}_{k}\right.  }\left(
Y\right)  \right]  \cdots\right]  ,E_{X_{2}\left\vert \mathbf{X}_{1}\right.
}\left[  \cdots\left[  E_{\mathbf{Y}\left\vert \mathbf{X}_{k}\right.  }\left(
Y\right)  \right]  \cdots\right]  \right) \\
&  +E_{X_{1}}\left[  Cov_{X_{2}\left\vert \mathbf{X}_{1}\right.  }\left(
E_{X_{3}\left\vert \mathbf{X}_{2}\right.  }\left[  \cdots\left[
E_{\mathbf{Y}\left\vert \mathbf{X}_{k}\right.  }\left(  Y\right)  \right]
\cdots\right]  ,E_{X_{3}\left\vert \mathbf{X}_{2}\right.  }\left[
\cdots\left[  E_{\mathbf{Y}\left\vert \mathbf{X}_{k}\right.  }\left(
Y\right)  \right]  \cdots\right]  \right)  \right] \\
&  \vdots\\
&  +E_{X_{1}}\left[  E_{X_{2}\left\vert \mathbf{X}_{1}\right.  }\left[
E_{X_{3}\left\vert \mathbf{X}_{2}\right.  }\left[  \cdots\left[
Cov_{\mathbf{Y}\left\vert \mathbf{X}_{k}\right.  }\left(  Y,Y\right)  \right]
\cdots\right]  \right]  \right] \\
&  =V_{X_{1}}\left(  E_{X_{2}\left\vert \mathbf{X}_{1}\right.  }\left[
\cdots\left[  E_{\mathbf{Y}\left\vert \mathbf{X}_{k}\right.  }\left(
Y\right)  \right]  \cdots\right]  ,E_{X_{2}\left\vert \mathbf{X}_{1}\right.
}\left[  \cdots\left[  E_{\mathbf{Y}\left\vert \mathbf{X}_{k}\right.  }\left(
Y\right)  \right]  \cdots\right]  \right) \\
&  +E_{X_{1}}\left[  V_{X_{2}\left\vert \mathbf{X}_{1}\right.  }\left(
E_{X_{3}\left\vert \mathbf{X}_{2}\right.  }\left[  \cdots\left[
E_{\mathbf{Y}\left\vert \mathbf{X}_{k}\right.  }\left(  Y\right)  \right]
\cdots\right]  ,E_{X_{3}\left\vert \mathbf{X}_{2}\right.  }\left[
\cdots\left[  E_{\mathbf{Y}\left\vert \mathbf{X}_{k}\right.  }\left(
Y\right)  \right]  \cdots\right]  \right)  \right] \\
&  \vdots\\
&  +E_{X_{1}}\left[  E_{X_{2}\left\vert \mathbf{X}_{1}\right.  }\left[
E_{X_{3}\left\vert \mathbf{X}_{2}\right.  }\left[  \cdots\left[
V_{\mathbf{Y}\left\vert \mathbf{X}_{k}\right.  }\left(  Y,Y\right)  \right]
\cdots\right]  \right]  \right]  \text{.}%
\end{align*}
This is the generalized law of total variance.$\blacksquare$

\section{Conclusion}

A generalization of the law of total covariance was given. It was proven that
the generalized law of total variance is a special case of the generalized law
of the total covariance. Examples were referenced.

\end{document}